\documentclass[leqno, 11pt]{article}
\usepackage{amsfonts}
\usepackage{graphics,graphicx,amssymb}
\usepackage{amsmath,amsthm}
\numberwithin{equation}{section}

\newtheorem{corollary}{Corollary}[section]
\newtheorem{proposition}{Proposition}[section]

\newtheorem{definition}{Definition}[section]
\newcommand{\tr}{{\rm tr}}

\newfont{\got}{eufm9 scaled 1095}

\newfont{\w}{msbm9 scaled\magstep1}

\newcommand{\Div}{\mathrm{div}}

\begin{document}

\title{On a Class Almost Contact Manifolds\\ with Norden Metric}

\footnotetext[1]{This work is partially supported by The Fund for
Scientific Research of the University of Plovdiv, Bulgaria, Project
RS09-FMI-003.}

\author{Marta Teofilova}

\date{}

\maketitle

\begin{abstract}
Certain curvature properties and scalar invariants of the manifolds
belonging to one of the main classes almost contact manifolds with
Norden metric are considered. An example illustrating the obtained
results is given and studied.

\end{abstract}

\section*{Introduction}

The geometry of the almost contact manifolds with Norden metric
($B$-metric) is a natural extension of the geometry of the almost
complex manifolds with Norden metric ($B$-metric) in the odd
dimensional case.

Almost contact manifolds with Norden metric are introduced in
\cite{Ga-Mi}. Eleven basic classes of these manifolds are
characterized there according to the properties of the covariant
derivatives of the almost contact structure.

In this work we focus our attention on one of the basic classes
almost contact manifolds with Norden metric, namely the class
$\mathcal{F}_{11}$. We study some curvature properties and relations
between certain scalar invariants of the manifolds belonging to this
class. In the last section we illustrate the obtained results by
constructing and studying an example of an
$\mathcal{F}_{11}$-manifold on a Lie group.

\section{Preliminaries}

Let $M$ be a $(2n+1)$-dimensional smooth manifold, and let
$(\varphi, \xi, \eta)$ be an almost contact structure on $M$, i.e.
$\varphi$ is an endomorphism of the tangent bundle of $M$, $\xi$ is
a vector field, and $\eta$ is its dual 1-form such that
\begin{equation}\label{1-1}
\varphi^2 = -\mathrm{Id} + \eta \otimes\xi,\qquad \eta(\xi)=1.
\end{equation}
Then, $(M,\varphi,\xi,\eta)$ is called an \emph{almost contact
manifold}.

We equip $(M,\varphi,\xi,\eta)$ with a compatible pseudo-Riemannian
metric $g$ satisfying
\begin{equation}\label{1-2}
g(\varphi x,\varphi y) = -g(x,y) + \eta(x)\eta(y)
\end{equation}
for arbitrary $x$, $y$ in the Lie algebra $\mathfrak{X}(M)$ of the
smooth vector fields on $M$. Then, $g$ is called a \emph{Norden
metric} ($B$-\emph{metric}), and $(M,\varphi,\xi,\eta,g)$ is called
an \emph{almost contact manifold with Norden metric}.

From (\ref{1-1}), $(\ref{1-2})$ it follows $\varphi \xi=0$, $\eta
\circ \varphi =0$, $\eta(x)=g(x,\xi)$, $g(\varphi x,y)=g(x,\varphi
y)$.

The associated metric $\tilde{g}$ of $g$ is defined by
$\tilde{g}(x,y)=g(x,\varphi y) + \eta(x)\eta(y)$ and is a Norden
metric, too. Both metrics are necessarily of signature $(n+1,n)$.

Further, $x,y,z,u$ will stand for arbitrary vector fields in
$\mathfrak{X}(M)$.

Let $\nabla$ be the Levi-Civita connection of $g$. The fundamental
tensor $F$ of type (0,3) is defined by
\begin{equation}\label{F}
F(x,y,z) = g\left((\nabla_x \varphi )y,z\right)
\end{equation}
and has the properties
\begin{equation}\label{Fp}
\begin{array}{l}
F(x,y,z)=F(x,z,y),\medskip\\
F(x,\varphi y,\varphi z)=F(x,y,z) -
F(x,\xi,z)\eta(y)-F(x,y,\xi)\eta(z).
\end{array}
\end{equation}
From the last equation and $\varphi \xi=0$ it follows
$F(x,\xi,\xi)=0$.

Let $\{e_i,\xi\}$ ($i=1,2,...,2n$) be a basis of the tangent space
$T_pM$ at an arbitrary point $p$ of $M$, and $g^{ij}$ be the
components of the inverse matrix of $(g_{ij})$ with respect to
$\{e_i,\xi\}$. The following 1-forms are associated with $F$:
\begin{equation}\label{teta}
\begin{array}{ll}
\theta(x)=g^{ij}F(e_i,e_j,x),\qquad&
\theta^{\ast}(x)=g^{ij}F(e_i,\varphi e_j,x),\medskip\\
\omega(x)=F(\xi,\xi,x),\qquad &\omega^\ast=\omega\circ\varphi.
\end{array}
\end{equation}
We denote by $\Omega$ the vector field corresponding to $\omega$,
i.e. $\omega(x)=g(x,\Omega)$.

The Nijenhuis tensor $N$ of the almost contact structure
$(\varphi,\xi,\eta)$ is defined by \cite{Sa}
$N(x,y)=[\varphi,\varphi](x,y)+d\eta(x,y)\xi$, i.e.
\begin{equation*}
N(x,y)=\varphi^2[x,y] + [\varphi x,\varphi y]- \varphi[\varphi x,y]-
\varphi [x,\varphi y] + (\nabla_x \eta)y.\xi - (\nabla_y \eta)x.\xi
\end{equation*}

In terms of the covariant derivatives of $\varphi$ and $\eta$ the
tensor $N$ is expressed as follows
\begin{equation}\label{Nfi}
\begin{array}{l}
N(x,y) = (\nabla_{\varphi x}\varphi) y - (\nabla_{\varphi y}\varphi
)x - \varphi(\nabla_x \varphi)y + \varphi(\nabla_y \varphi)x
\medskip\\
\phantom{N(x,y)} + (\nabla_x \eta) y. \xi - (\nabla_y \eta) x. \xi,
\end{array}
\end{equation}
where $(\nabla_x \eta)y=F(x,\varphi y,\xi)$.

The almost contact structure in said to be integrable if $N=0$. In
this case the almost contact manifold is called \emph{normal}
\cite{Sa}.

A classification of the almost contact manifolds with Norden metric
is introduced in \cite{Ga-Mi}. This classification consists of
eleven basic classes $\mathcal{F}_i$ ($i=1,2,...,11$) characterized
according to the properties of $F$. The special class
$\mathcal{F}_0$ of the $\varphi$-K\"ahler-type almost contact
manifolds with Norden metric is given by the condition $F=0$
($\nabla\varphi=\nabla\xi=\nabla\eta=0$). The classes for which $F$
is expressed explicitly by the other structural tensors are called
\emph{main classes}.

In the present work we focus our attention on one of the main
classes of these manifolds, namely the class $\mathcal{F}_{11}$,
which is defined by the characteristic condition \cite{Ga-Mi}
\begin{equation}\label{F11}
F(x,y,z)=\eta(x)\{\eta(y)\omega(z)+\eta(z)\omega(y)\}.
\end{equation}

By (\ref{teta}) and (\ref{F11}) we get that on a
$\mathcal{F}_{11}$-manifold $\theta=\omega$, $\theta^\ast=0$. We
also have
\begin{equation}\label{ho}
(\nabla_x\omega^\ast)y = (\nabla_x \omega)\varphi y +
\eta(x)\eta(y)\omega(\Omega).
\end{equation}

The 1-forms $\omega$ and $\omega^\ast$ are said to be closed if
$\text{d}\omega=\text{d}\omega^\ast=0$. Since $\nabla$ is symmetric,
necessary and sufficient conditions for $\omega$ and $\omega^\ast$
to be closed are
\begin{equation}\label{omega}
(\nabla_x\omega)y=(\nabla_y\omega)x,\qquad(\nabla_x\omega)\varphi
y=(\nabla_y\omega)\varphi x.
\end{equation}

The curvature tensor $R$ of $\nabla$ is defined as usually by
\begin{equation}\label{R}
R(x,y)z=\nabla_x\nabla_y z - \nabla_y\nabla_x z - \nabla_{[x,y]}z,
\end{equation}
and its corresponding tensor of type (0,4) is given by
$R(x,y,z,u)=g(R(x,y)z,u)$. The Ricci tensor $\rho$ and the scalar
curvatures $\tau$ and $\tau^{\ast}$ are defined by, respectively
\begin{equation}
\rho(y,z)=g^{ij}R(e_{i},y,z,e_{j}),\qquad
\tau=g^{ij}\rho(e_{i},e_{j}),\qquad \tau^{\ast}=g^{ij}\rho
(e_{i},\varphi e_{j}). \label{Ricci-tao}
\end{equation}

The tensor $R$ is said to be of $\varphi$-K\"ahler-type if
\begin{equation}\label{KR}
R(x,y,\varphi z,\varphi u) = - R(x,y,z,u).
\end{equation}

Let $\alpha=\{x,y\}$ be a non-degenerate 2-section spanned by the
vectors $x,y\in T_pM$, $p\in M$. The sectional curvature of $\alpha$
is defined by
\begin{equation}\label{sec}
k(\alpha;p)=\frac{R(x,y,y,x)}{\pi_1(x,y,y,x)},
\end{equation}
where $\pi_1(x,y,z,u) = g(y,z)g(x,u) - g(x,z)g(y,u)$.

In \cite{Na-Gri} there are introduced the following special sections
in $T_pM$: a $\xi$-section if $\alpha=\{x,\xi\}$, a
$\varphi$-holomorphic section if $\varphi \alpha = \alpha$ and a
totally real section if $\varphi \alpha \perp \alpha$ with respect
to $g$.

The square norms of $\nabla \varphi$, $\nabla \eta$ and $\nabla\xi$
are defined by, respectively \cite{Ko-No}:
\begin{equation}\label{n}
\begin{array}{l}
||\nabla\varphi||^2 =
g^{ij}g^{ks}g\left((\nabla_{e_i}\varphi)e_{k},(\nabla_{e_j}\varphi)e_{s}\right),\medskip\\
||\nabla \eta||^2=||\nabla\xi||^2 = g^{ij}g^{ks}(\nabla_{e_i}\eta)
e_{k}(\nabla_{e_j}\eta)e_{s}.
\end{array}
\end{equation}

We introduce the notion of an isotropic K\"ahler-type almost contact
manifold with Norden metric analogously to \cite{Gri-Man-Mek2}.

\begin{definition}\label{def1} \emph{An almost contact manifold with Norden metric is
called} isotropic K\"ahlerian \emph{if
$||\nabla\varphi||^2=||\nabla\eta||^2=0$ (
$||\nabla\varphi||^2=||\nabla\xi||^2=0$).}
\end{definition}

\section{Curvature properties of $\mathcal{F}_{11}$-manifolds}

In this section we obtain relations between certain scalar
invariants on $\mathcal{F}_{11}$-manifolds with Norden metric and
give necessary and sufficient conditions for such manifolds to be
isotropic K\"aherian.

First, by help of (\ref{F}), (\ref{Fp}), (\ref{Nfi}), (\ref{F11}),
(\ref{n}) and direct computation we obtain
\begin{proposition}\label{14}
On a $\mathcal{F}_{11}$-manifold it is valid
\begin{equation}
||\nabla\varphi||^2=-||N||^2=-2||\nabla\eta||^2=2\omega(\Omega).
\end{equation}
\end{proposition}
Then, (\ref{14}) and Definition \ref{def1} yield
\begin{corollary}\label{cor1}
On a $\mathcal{F}_{11}$-manifold the following conditions are
equivalent\emph{:}
\begin{itemize}
\item[\emph{(i)}] the manifold is isotropic
K\"ahlerian\emph{;}
\item[\emph{(ii)}] the vector $\Omega$ is isotopic, i.e. $\omega(\Omega)=0$\emph{;}
\item[\emph{(iii)}] the Nijenhuis tensor $N$ is isotropic.
\end{itemize}
\end{corollary}

It is known that the almost contact structure satisfies the Ricci
identity, i.e.
\begin{equation}\label{Ricci}
\begin{array}{l}
(\nabla_x\nabla_y\varphi)z - (\nabla_y\nabla_x\varphi)z =
R(x,y)\varphi z -\varphi R(x,y)z,\medskip\\ (\nabla_x\nabla_y\eta)z
- (\nabla_y\nabla_x\eta)z = -\eta(R(x,y)z).
\end{array}
\end{equation}
Then, taking into account the definitions of $\varphi$, $F$, and
$\nabla g=0$, the equalities (\ref{Ricci}) imply
\begin{equation}\label{h1}
\begin{array}{c}
(\nabla_xF)(y,z,\varphi u) - (\nabla_yF)(x,z,\varphi u) =
R(x,y,z,u)\medskip\\+ R(x,y,\varphi z,\varphi u)-
R(x,y,z,\xi)\eta(u),
\end{array}
\end{equation}
\begin{equation}\label{h2}
(\nabla_xF)(y,\varphi z,\xi) - (\nabla_yF)(x,\varphi z,\xi) =
-R(x,y,z,\xi).
\end{equation}

By (\ref{F}), (\ref{h1}) and (\ref{h2}) we get
\begin{equation}\label{RR}
R(x,y,\varphi z,\varphi u) = - R(x,y,z,u) + \psi_4(S)(x,y,z,u),
\end{equation}
where the tensor $\psi_4(S)$ is defined by \cite{Ma-4}
\begin{equation}\label{psi}
\begin{array}{l}
\psi_4(S)(x,y,z,u) = \eta(y)\eta(z)S(x,u) -
\eta(x)\eta(z)S(y,u)\medskip\\
\phantom{\psi_4(S)(x,y,z,u)} + \eta(x)\eta(u)S(y,z) -
\eta(y)\eta(u)S(x,z).
\end{array}
\end{equation}
and
\begin{equation}\label{S}
S(x,y)=(\nabla_x\omega)\varphi y - \omega(\varphi x)\omega(\varphi
y).
\end{equation}

Then, the following holds
\begin{proposition}
On a $\mathcal{F}_{11}$-manifold we have
\begin{equation*}
\tau+\tau^{\ast\ast} =2\Div(\varphi\Omega)=2\rho(\xi,\xi),
\end{equation*}
where $\tau^{\ast\ast}=g^{is}g^{jk}R(e_i,e_j,\varphi e_k,\varphi
e_s)$.
\begin{proof}
The truthfulness of the statement follows from (\ref{Ricci-tao}) and
(\ref{RR}) by straightforward computation.
\end{proof}
\end{proposition}

Having in mind (\ref{KR}) and (\ref{RR}), we conclude that the
curvature tensor on a $\mathcal{F}_{11}$-manifold is of
$\varphi$-K\"aher-type if and only if $\psi_4(S)=0$. Because of
(\ref{psi}) the last condition holds true iff $S=0$. Then, taking
into account (\ref{ho}) and (\ref{S}) we prove
\begin{proposition}
The curvature tensor of a $\mathcal{F}_{11}$-manifold with Norden
metric is $\varphi$-K\"aherian iff
\begin{equation}\label{S0}
(\nabla_x\omega^\ast)y = \eta(x)\eta(y)\omega(\Omega)+ \omega^\ast(
x)\omega^\ast(y).
\end{equation}
\end{proposition}
The condition (\ref{S0}) implies $\text{d}\omega^\ast=0$, i.e.
$\omega^\ast$ is closed.

\section{An example}

In this section we present and study a $(2n+1)$-dimensional example
of a $\mathcal{F}_{11}$-manifold constructed on a Lie group.

Let $G$ be a $(2n+1)$-dimensional real connected Lie group, and
$\mathfrak{g}$ be its corresponding Lie algebra. If
$\{x_0,x_1,...,x_{2n}\}$ is a basis of left-invariant vector fields
on $G$, we define a left-invariant almost contact structure
$(\varphi,\xi,\eta)$ by
\begin{equation}\label{fi}
\begin{array}{llll}
\varphi x_{i}= x_{i+n},\quad& \varphi x_{i+n}=-x_i,\quad& \varphi
x_{0}=0,\quad& i=1,2,...,n,\medskip\\
\xi = x_{0},\quad & \eta(x_0)=1, \quad& \eta(x_{j})=0, \quad&
j=1,2,...,2n.
\end{array}
\end{equation}
We also define a left-invariant pseudo-Riemannian metric $g$ on $G$
by
\begin{equation}\label{g}
\begin{array}{l}
g(x_0,x_0)=g(x_i,x_i)=-g(x_{i+n},x_{i+n})=1,\quad i=1,2,...,n,\medskip\\
g(x_j,x_k)=0,\quad j\neq k,\quad j,k=0,1,,...,2n.
\end{array}
\end{equation}
Then, according to (\ref{1-1}) and (\ref{1-2}),
$(G,\varphi,\xi,\eta,g)$ is an almost contact manifold with Norden
metric.

Let the Lie algebra $\mathfrak{g}$ of $G$ be given by the following
non-zero commutators
\begin{equation}\label{Lie}
[x_i,x_0]=\lambda_i x_0,\qquad i=1,2,...,2n,
\end{equation}
where $\lambda_i\in\mathbb{R}$. Equalities (\ref{Lie}) determine a
$2n$-parametric family of solvable Lie algebras.

Further, we study the manifold $(G,\varphi,\xi,\eta,g)$ with Lie
algebra $\mathfrak{g}$ defined by (\ref{Lie}). The well-known
Koszul's formula for the Levi-Civita connection of $g$ on $G$, i.e.
the equality
\begin{equation*}
2g(\nabla_{x_i} x_j,x_k) = g([x_i,x_j],x_k) +
g([x_k,x_i],x_j)+g([x_k,x_j],x_i),
\end{equation*}
implies the following components of the Levi-Civita connection:
\begin{equation}
\begin{array}{l}\label{nabla}
\nabla_{x_{i}}x_{j}=\nabla_{x_{i}}\xi=0,\qquad
\nabla_{\xi}x_{i}=-\lambda_i\xi,\quad i,j=1,2,...,2n,\medskip\\
\nabla_{\xi}\xi = \sum_{k=1}^n (\lambda_k x_k -
\lambda_{k+n}x_{k+n}).
\end{array}
\end{equation}
Then, by (\ref{F}) and (\ref{nabla}) we obtain the essential
non-zero components of $F$:
\begin{equation}\label{Fi}
F(\xi,\xi,x_i)=\omega(x_i)= -\lambda_{i+n},\qquad
F(\xi,\xi,x_{i+n})=\omega(x_{i+n})=\lambda_i,
\end{equation}
for $i=1,2,...,n$. Hence, by (\ref{F11}) and (\ref{Fi}) we have
\begin{proposition}
The almost contact manifold with Norden metric
$(G,\varphi,\xi,\eta,g)$ defined by (\ref{fi}), (\ref{g}) and
(\ref{Lie}) belongs to the class $\mathcal{F}_{11}$.
\end{proposition}

Moreover, by (\ref{omega}), (\ref{nabla}) and (\ref{Fi}) we
establish that the considered manifold has closed 1-forms $\omega$
and $\omega^\ast$.

Taking into account (\ref{nabla}) and (\ref{R}) we obtain the
essential non-zero components of the curvature tensor as follows
\begin{equation}\label{Rijks}
R(x_i,\xi,\xi,x_j) = -\lambda_i\lambda_j,\quad i,j=1,2,...,2n.
\end{equation}
By (\ref{Rijks}) it follows that $R(x_i,x_j,\varphi x_k,\varphi
x_s)=0$ for all $i,j,k,s=0,1,...,2n$. Then, according to (\ref{RR})
and (\ref{S}) we get
\begin{proposition}
The curvature tensor and the Ricci tensor of the
$\mathcal{F}_{11}$-manifold $(G,\varphi,\xi,\eta,g)$ defined by
(\ref{fi}), (\ref{g}) and (\ref{Lie}) have the form, respectively
\begin{equation*}
R = \psi_4(S),\qquad \rho(x,y) = \eta(x)\eta(y) \tr S + S(x,y),
\end{equation*}
where $S$ is defined by (\ref{S}) and $\tr S = \Div (\varphi
\Omega)$.
\end{proposition}

We compute the essential non-zero components of the Ricci tensor as
follows
\begin{equation}
\begin{array}{l}\label{Rij}
\rho(x_i,x_j)=-\lambda_i\lambda_j,\quad i=1,2,...,2n,\medskip\\
\rho(\xi,\xi) = -\sum_{k=1}^n\left(\lambda_k^2 -
\lambda_{k+n}^2\right).
\end{array}
\end{equation}
By (\ref{Ricci-tao}) and (\ref{Rij}) we obtain the curvatures of the
considered manifold
\begin{equation}
\tau = -2\sum_{k=1}^n\left(\lambda_k^2 -
\lambda_{k+n}^2\right),\qquad \tau^\ast=-2\sum_{k=1}^n
\lambda_k\lambda_{k+n}.
\end{equation}

Let us consider the characteristic 2-sections $\alpha_{ij}$ spanned
by the vectors $\{x_i,x_j\}$: $\xi$-sections $\alpha_{0,i}$
($i=1,2,...,2n$), $\varphi$-holomorphic sections $\alpha_{i,i+n}$
($i=1,2,...,n$), and the rest are totally real sections. Then, by
(\ref{sec}), (\ref{fi}) and (\ref{Rijks}) it follows
\begin{proposition}
The $\mathcal{F}_{11}$-manifold with Norden metric
$(G,\varphi,\xi,\eta,g)$ defined by (\ref{fi}), (\ref{g}) and
(\ref{Lie}) has zero totally real and $\varphi$-holomorphic
sectional curvatures, and its $\xi$-sectional curvatures are given
by
\begin{equation*}
k(\alpha_{0,i}) = -\frac{\lambda_i^2}{g(x_i,x_i)},\qquad i
=1,2,...,2n.
\end{equation*}
\end{proposition}

By (\ref{g}) and (\ref{Fi}) we obtain the corresponding vector
$\Omega$ to $\omega$ and its square norm
\begin{equation}\label{omega1}
\Omega = -\sum_{k=1}^n\left(\lambda_{k+n}x_k+\lambda_k
x_{k+n}\right),\qquad \omega(\Omega) =
-\sum_{k=1}^n\left(\lambda_k^2 - \lambda_{k+n}^2\right).
\end{equation}
Then, by (\ref{omega1}) and Corollary \ref{cor1} we prove
\begin{proposition}
The $\mathcal{F}_{11}$-manifold with Norden metric
$(G,\varphi,\xi,\eta,g)$ defined by (\ref{fi}), (\ref{g}) and
(\ref{Lie}) is isotropic K\"ahlerian iff the condition
$\sum_{k=1}^n\left(\lambda_k^2 - \lambda_{k+n}^2\right)=0$ holds.
\end{proposition}


\smallskip

\begin{tabbing}
  \emph{Marta Teofilova}\\
  \emph{University of Plovdiv}\\
  \emph{Faculty of Mathematics and Informatics}\\
  \emph{236 Bulgaria Blvd.}\\
  \emph{4003 Plovdiv, Bulgaria}\\
  \texttt{e-mail:\ marta@uni-plovdiv.bg}\\
  \end{tabbing}

\end{document}